\documentclass[12pt]{amsart}
\usepackage{mathrsfs}
\usepackage{amsfonts}
\usepackage{amsthm}
\usepackage{amssymb}
\usepackage{latexsym}
\newtheorem{theorem}{Theorem}

\newtheorem{corollary}{Corollary}

\newcommand{\dbar}{\overline{\partial}}
\newcommand{\dbars}{\overline{\partial}^*}

\begin{document}

\title[Complex tangential flows and compactness]{Complex tangential flows and compactness of the $\overline{\partial}$-Neumann operator}

\author{Samangi Munasinghe \and Emil J. Straube}
\address{Department of Mathematics\\
Texas A\&M University\\
College Station, Texas, 77843--3368}
\email{munasing@math.tamu.edu, straube@math.tamu.edu}

\thanks{2000 \emph{Mathematics Subject Classification}: 32W05}
\keywords{$\overline{\partial}$-Neumann operator, compactness,
complex tangential flows}
\thanks{Research supported in part by NSF grant DMS 0500842}

\date{July 25, 2006}

\begin{abstract}
We provide geometric conditions on the set of boundary points of infinite type of a smooth bounded pseudoconvex domain in $\mathbb{C}^{n}$ which imply that the $\overline{\partial}$-Neumann  operator is compact. These conditions are formulated in terms of certain short time flows in suitable complex tangential directions. It is noteworthy that compactness is \emph{not} established via the known potential theoretic sufficient conditions. Our results generalize to $\mathbb{C}^{n}$ the $\mathbb C^{2}$ results from \cite{St2}. 
\end{abstract}

\maketitle

\section{Introduction}

In \cite{St2}, the second author provided geometric sufficient conditions for compactness of the $\overline{\partial}$-Neumann operator on the boundary points of infinite type of a bounded smooth pseudoconvex domain in $\mathbb{C}^{2}$. In this paper, we study the situation in higher dimensions and obtain suitable generalizations of the results in \cite {St2}.

Let $\Omega$ be a bounded pseudoconvex domain in $\mathbb{C}^{n}$. The $\overline{\partial}$-Neumann operator $N_{q}$ on $(0,q)$-forms is the inverse of the complex Laplacian $\overline{\partial}\overline{\partial}^{*} + \overline{\partial}^{*}\overline{\partial}$ associated to the Dolbeault complex. For detailed information on the $\overline{\partial}$-Neumann problem and related questions, see e.g.\cite{FK,CS,LM,BS,McN2,St3}; compactness of $N_{q}$ is discussed in \cite{FS1, St3}.

Whether or not $N_{q}$ is compact is relevant in a number of situations. These include global regularity \cite{KN}, the Fredholm theory of Toeplitz operators \cite{HI}, and the (non)existence of solution operators to $\overline{\partial}$ with well-behaved solution kernels \cite{HL}. There are also interesting connections to the theory of Schr\"{o}dinger operators \cite{FS3,CF}. Catlin gave a sufficient condition, which he called property($P$), in \cite{C1}: there should exist, near the boundary, plurisubharmonic functions bounded between $0$ and $1$, with arbitrarily large Hessians. (The smoothness assumptions on the boundary of the domain were considerably weakened in \cite{HI}, and were completely removed in \cite{St}.) Property($P$) was studied in detail (under the name $B$-regularity) by Sibony in \cite{S} (see also \cite{S2}). On sufficiently regular domains, property($P$) is equivalent to a quantitative version of Oka's lemma (\cite{H}). A version of property($P$), called condition ($\widetilde{P}$), was introduced, and shown to still imply compactness, by McNeal in \cite{McN}. The uniform bound on the family of functions is replaced by a uniform bound on the gradient, measured in the metric induced by the complex Hessian of the functions. (Both $(P)$ and $(\widetilde{P})$ can also be formulated naturally at the level of $(0,q)$-forms; then $(P_{q}) \Rightarrow (P_{q+1})$, $(\widetilde{P}_{q}) \Rightarrow (\widetilde{P}_{q+1})$, and $(P_{q}) \Rightarrow (\widetilde{P}_{q})$, $1 \leq q \leq n$, see \cite{FS1,McN,St3}.) A sufficient condition that is intermediate, in a sense one can make precise (see the discussion in \cite{St3}), had appeared earlier in \cite{T}. 

On locally convexifiable domains, ($P_{q}$) and ($\widetilde{P}_{q}$) are equivalent, and equivalent to compactness of $N_{q}$, for $1 \leq q \leq n$. Moreover, the three properties are equivalent to a simple geometric condition, the absence of (germs of) $q$-dimensional varieties from the boundary. For this, see \cite{FS2,FS1}. Thus the potential theory, the analysis, and the geometry mesh perfectly on locally convexifiable domains. It is also known that on smooth bounded Hartogs domains in $\mathbb{C}^{2}$, compactness of $N_{1}$ is equivalent to ($P_{1}$) and to ($\widetilde{P}_{1}$) (\cite{CF,FS3}). However, it is well understood that the boundaries of convex (hence of locally convexifiable) domains do not exhibit some of the more intriguing aspects of the interaction with the ambient space that occur on general pseudoconvex boundaries (\cite{DA}). For example, matters concerning orders of contact are always decided by orders of contact of \emph{manifolds} (affine manifolds in the convex case) (\cite{McN1,BS1,Y,FS2}). A similar caveat applies in the case of domains in $\mathbb{C}^{2}$ (\cite{DA}). As a result, these facts give no clear indication of how much (or how little) room there is, in the general pseudoconvex case, between ($P$)/($\widetilde{P}$) and compactness of the $\overline{\partial}$-Neumann operator. (As far as we know, the exact relationship between ($P$) and ($\widetilde{P}$) is also unknown, but we do not address this question here.)

From the point of view of obstructions, the situation is as follows. $q$-dimensional varieties in the boundary are an obstruction for both property($P_{q}$) and condition($\widetilde{P}_{q}$) (\cite{S,FS1,St4} for ($P_{q}$), \cite{St3,St4} for ($\widetilde{P}_{q}$)). With respect to compactness of $N_{q}$, less is known. A $q$-dimensional complex manifold $M$ in the boundary of a smooth bounded pseudoconvex domain is known to be an obstruction to compactness of $N_{q}$, provided $M$ contains a point at which the domain is strictly pseudoconvex in the directions transverse to $M$ (\cite{Sa,SaSt}). It is open whether the conclusion holds without assuming that there is such a point. One would expect that a flatter boundary is even more favorable for noncompactness, but the methods of \cite{SaSt} do not seem strong enough to yield noncompactness without some additional assumption on how $M$ sits inside the boundary. On the other hand, both for ($P$)/($\widetilde{P}$) and for compactness of the $\overline{\partial}$-Neumann operator, there are obstructions more subtle than varieties in the boundary (\cite{S,M,FS1}). For a connection of some of these issues with properties of the Kobayashi metric, see \cite{KIM}.

\cite{St2} provided, for the first time, a method to prove compactness of the $\overline{\partial}$-Neumann operator that does not proceed via verifying property($P$) or condition ($\widetilde{P}$). That the dimension is two was only used in the application of so-called maximal estimates. Consequently, the results of \cite{St2} hold more generally on domains in $\mathbb{C}^{n}$, $n \geq 2$, where such estimates hold, or, equivalently, on domains where all the eigenvalues of the Levi form are comparable {\cite{D1,D2}. However, for the problem of compactness of the $\overline{\partial}$-Neumann operator, this assumption is too restrictive. It excludes, for example, the situation where the Levi form has at most one degenerate eigenvalue (see Remark 5 below). But this assumption on the Levi form has been shown to be a useful generalization of the case of $\mathbb{C}^{2}$ in the context of compactness (\cite{Sa,SaSt}).

The obvious examples that satisfy the assumptions in Theorem \ref{main} below also satisfy property($\widetilde{P}$); we do not know whether the theorem can actually furnish examples of domains where the $\overline{\partial}$-Neumann operator is compact, but where $(\widetilde{P}$) fails. But just as in \cite{St2}, we obtain a simple geometric proof of compactness in these cases. Moreover,the assumptions are in some instances `minimal': they are necessary modulo the size of certain balls; see Remark 6 below for details.

We will only consider the case $q=1$ in the remainder of this paper. This is the main case in terms of understanding compactness. But note that compactness of $N_{1}$ implies compactness of $N_{q}$ for $q >1$; this is an observation due to McNeal (\cite{McN3}, see also the proof of Lemma 2 in \cite{St5} for a related argument).

The bulk of this paper represents a portion of the first author's Ph.D. thesis (\cite{MU}) written at Texas A\&M University under the direction of the second author.  

\section{Results}

If $Z$ is a (real) vector field defined in some open subset of the boundary (or of $\mathbb{C}^{n}$), we denote by $\mathcal{F}_{Z}^{t}$ the flow generated by $Z$. We use the notion of finite or infinite type of D'Angelo \cite{DA}. For a boundary point $\zeta$, we denote by $\lambda_{0}(\zeta)$ the smallest eigenvalue of the Levi form of the boundary at $\zeta$. Since the domains in question are pseudoconvex, $\lambda_{0}(\zeta) \geq 0$. For a set of real vector fields $T_{1}, \cdots, T_{m}$ on an open subset of $b\Omega$ (or of $\mathbb{C}^{n}$), we define $span_{\mathbb{R}}(T_{1}, \cdots, T_{m})$ to be the set of all linear combinations of   $T_{1}, \cdots, T_{m}$ whose coefficients are real valued functions (not necessarily constants).

\begin{theorem}\label{main}
Let $\Omega$ be a $C^\infty$-smooth bounded pseudoconvex domain in
$\mathbb{C}^n$.  Denote by $K$ the set of boundary points of infinite type.  
Assume that there exist smooth complex tangential vector fields
$X_1, \dots X_m$, defined on b$\Omega$ near $K$, so that
$H_{\rho}(X_i(\zeta), \overline{X_i(\zeta)})$ $\leq C{\lambda}_0(\zeta)$,
for some constant $C$, a sequence 
${\lbrace {\epsilon_j} \rbrace}_{j=1}^{\infty}$ with 
$\underset{j \rightarrow \infty}{\text{lim}}\epsilon_j=0$, 
and constants $C_1$, $C_2 >0$, $C_3$ with $1 \leq C_3 < \frac{n+1}{n}$, so that the following holds.
For every $j \in \mathbb{N}$ and $p \in K$ there is a real 
vector field $Z_{p,j} \in \text{span}_{\mathbb{R}}(\text{Re}X_1, \text{Im}X_1,
\dots , \text{Re}X_m, \text{Im}X_m)$ of unit length, defined in some neighborhood of $p$ in b$\Omega$ with ${\text{max}} \lvert \text{div} Z_{p,j} \rvert \leq C_1 $, such that
$\mathcal{F}_{Z_{p}}^{\epsilon_j}{\left(B(p, C_2(\epsilon_j)^{C_3})\cap K\right)} \subseteq b\Omega \setminus K$.  
Then the $\bar{\partial}$-Neumann operator $N_{1}$ on $\Omega$ is compact.
\end{theorem}
\smallskip
\emph{Remark 1}:
In $\mathbb{C}^2$, take $m=1$ and $X_1 = L$, where $L$ is a smooth nonvanishing complex tangential vector field of of type $(1,0)$ on b$\Omega$. Then the condition on the vector fields $Z_{p,j}$ becomes simply that they be complex tangential, as in \cite{St2}. Theorem \ref{main} therefore generalizes the main result in \cite{St2}.

\medskip
\emph{Remark 2}:
In comparison to the main result in \cite{St2}, Theorem \ref{main} contains the additional assumption that there is a family of vector fields
$X_1 \dots X_m$ satisfying
\begin{equation}\label{new-cond}
H_{\rho}(X_j(\zeta), \overline{X_j(\zeta)}) \leq C{\lambda}_0(\zeta) \; ,
\end{equation}   
such that the vector fields $Z_{p,j}$ are contained in the linear span (over $\mathbb{R}$, in the sense of the definition above) of the real and imaginary parts of these fields (as opposed to just being complex tangential). Without some assumption on the fields $Z_{p,j}$ more restrictive than being complex tangential, such as the one made here, the result does not generalize to $\mathbb{C}^n$. To see this, consider a smooth bounded convex domain in $\mathbb{C}^3$ which is strictly convex, except for an analytic (affine) disc in the boundary.  Then one can flow along complex tangential directions from points of the disc into the set of strictly (pseudo) convex boundary points as required in the second part of the assumption in the theorem.  Nonetheless, because there is a disc in the boundary, the $\dbar$-Neumann operator on $(0,1)$-forms is not compact on such a domain (see \cite{FS2}).

\medskip
\emph{Remark 3}:
In the example in Remark 2, the value of the Levi form in the direction parallel to the disc apparently goes to zero faster, upon approach to the disc, than the value in the direction transverse to the disc. Thus (\ref{new-cond}) cannot hold for a vector field transverse to the disc, and Theorem \ref{main} does not apply. It would be very desirable to have a direct proof of this in general, i.e. a proof that the assumptions in Theorem \ref{main} exclude discs from the boundary (or, less likely, a counterexample). This is relevant because the question whether in general, discs in the boundary are obstructions to compactness is open (compare the discussion in the introduction for what is known). Note that in $\mathbb{C}^{2}$, such discs are known to be obstructions to compactness (\cite{FS1}); in this case, it is also obvious that the assumptions in the theorem exclude discs from the boundary.

\medskip
\emph{Remark 4}: Because compactness of the $\overline{\partial}$-Neumann operator follows from a compactness estimate for forms supported in fixed, but possibly small, neighborhoods of boundary points (see e.g. \cite{FS1}), there is a version of Theorem \ref{main} where the assumptions are localized. For this, see the discussion in Example 3 below.

\medskip
In Theorem \ref{main}, one would like to have a collection of vector fields $\{X_{j}\}_{j=1}^{m}$ such that at each point $p$ of $K$, $span_{\mathbb{R}}\{ReX_{1}(p), ImX_{1}(p), \cdots , ReX_{m}(p), ImX_{m}(p)\}$ is as `big as possible', thus putting the least restrictions on the fields $Z_{p,j}$. On the other hand, this needs to be balanced with the requirement (\ref{new-cond}). We now discuss some examples.

\smallskip
\emph{Example 1}:
When the eigenvalues of the Levi form are all comparable, any finite collection of complex tangential vector fields $X_1, \dots, X_{m}$ will satisfy condition (\ref{new-cond}). Taking a collection which at each point $p\in K$ spans all of $T^{\mathbb{C}}_{b\Omega}(p)$, we see that in this case, this part of the assumptions in Theorem \ref{main} reduces to $Z_{p,j}$ complex tangential (as in $\mathbb{C}^{2}$).
Domains where all the eigenvalues of the Levi form are comparable play a special role in the theory of the $\overline{\partial}$-Neumann problem; certain estimates, called `maximal estimates', hold. This class of domains was studied in detail in \cite{D1, D2}.

\smallskip
\emph{Example 2}:
Assume there exists a smooth complex tangential vector field $X_{1}$ near $K$ such that for $\zeta \in K$, $X_{1}(\zeta)$ is an eigenvector associated with the smallest eigenvalue of the Levi form at $\zeta$. This vector field trivially satisfies the condition (\ref{new-cond}). Then the assumption in Theorem \ref{main} requires that $Z_{p,j}(\zeta)$ is in the two real-dimensional plane spanned by $X_{1}(\zeta)$ for all $\zeta$.

\smallskip
\emph{Example 3}:
Assume that the Levi form has at most one degenerate eigenvalue at each point of $K$ (hence near $K$). Fix a point $p \in K$. Choose an $(n-2)$-dimensional subspace of $T^{\mathbb{C}}_{b\Omega}(p)$ on which the Levi form is strictly positive, and choose a basis. Extending the basis vectors to local sections of $T^{\mathbb{C}}_{b\Omega}$ gives vector fields (defined near $p$) $Y_{2}, \cdots , Y_{n-1}$. In a neighborhood of $p$, the Levi form is strictly positive on the span of $Y_{2}, \cdots , Y_{n-1}$. As a consequence, at each point there is a unique one dimensional subspace of $T^{\mathbb{C}}_{b\Omega}$ that is orthogonal to this span with respect to the Levi form. Indeed, if $Y_{1}$ is such that $Y_{1}, \cdots , Y_{n-1}$ is a basis for $T^{\mathbb{C}}_{b\Omega}$ (near $p$), then $X_{1} = Y_{1} + 
b_{2}Y_{2}  \cdots + b_{n-1}Y_{n-1}$, where $b_{j}=-H_{\rho}(Y_{1},\overline{Y_{j}})/H_{\rho}(Y_{j},\overline{Y_{j}})$, $j=2,  \cdots,  n-1$, will span this subspace. (This was observed already in \cite{MA}, Lemma 2.1.) Note that when $p \in K$, then the eigenvector of the Levi form at $p$ associated with the eigenvalue zero is orthogonal to all of $T^{\mathbb{C}}_{b\Omega}(p)$ with respect to the Levi form (because the Levi form is positive definite). Therefore,by uniqueness, $X_{1}(p)$ is an eigenvector of the Levi form at $p$ with eigenvalue zero. As a result, $H_{\rho}(X_{1}(\zeta),\overline{X_{1}(\zeta)}) \leq C\lambda_{0}(\zeta)$ for $\zeta$ close enough to $p$ (by continuity), that is, (\ref{new-cond}) holds for the family consisting of the single field $X_{1}$. We note that when $\zeta$ is a strictly pseudoconvex point, $X_{1}(\zeta)$ need not be an eigenvector of the Levi form. By multiplying with a cutoff function that is identically equal to one near $p$, we may assume that $X_{1}(\zeta)$ is defined on all of $b\Omega$, with (\ref{new-cond}) still valid. Of course, this is at the expense of having trivial span on a big set. We may proceed in two ways. We can cover $K$ with finitely many open sets $U_{1}, \cdots, U_{s}$ on which the cutoff functions multiplying these local fields are one, and then add finitely many of the fields. However, the resulting field may still vanish at some points of $K$. Alternatively, we can take advantage of the fact that compactness localizes: it suffices to prove a compactness estimate for forms supported in (small) neighborhoods of boundary points (as long as the neighborhood do not depend on the $\epsilon$ in the compactness estimate). The proof of Theorem \ref{main}, using the field $X_{j}$, gives a compactness estimate for forms whose support meets the boundary in one of the $U_{j}$'s. Since the $U_{j}$'s cover $K$, the result is a compactness estimate for all forms in $dom(\overline{\partial}) \cap dom(\overline{\partial}^{*})$ (since away from $K$, we have subelliptic estimates).

\smallskip
\emph{Remark 5}: On the domains discussed in Example 3, the Levi form is locally diagonalizable (use Gram-Schmitt to orthonormalize $Y_{2}, \cdots, Y_{n-1}$ with respect to the Levi form). Derridj showed in (\cite{D2}, Theorem 7.1), that if maximal estimates hold at $p \in b\Omega$, and $p$ is a weakly pseudoconvex point, then the Levi form of $\Omega$ cannot be diagonalizable near $p$ when $\Omega$ is a domain in $\mathbb{C}^n$ with $n \geq 3$.  Therefore, the examples (1) and (3) are mutually exclusive (when $n \geq 3$). 

\smallskip
\emph{Remark 6}: It is interesting to note that in Example 3, the assumptions in the theorem are `minimal', that is, they are necessary modulo the size of the balls $B(p,C_{2}(\epsilon_{j})^{C_{3}})$. The discussion is analogous to the one in Remarks 2 and 3 in \cite{St2}, but using recent results from \cite{SaSt}. For $p \in K$, let $X_{1}$ be the complex tangential field from Example 3 above (defined near $p$, with $H_{\rho}(X_{1}(p),\overline{X_{1}(p)})=0$). Denote by $T^{\theta}$ the field $T^{\theta}=cos(\theta)ReX_{1} + sin(\theta)ImX_{1}$. For $\zeta$ near $p$, set $M_{\zeta,\theta}=\{\mathcal{F}_{T^{\theta}}^{t}(\zeta)|0\leq \theta \leq 2\pi, 0\leq t \leq t_{0}\}$. Then $M_{\zeta,\theta}$ is a smooth two dimensional submanifold of the boundary. Because $N_{1}$ is compact, the boundary contains no analytic discs (since the Levi form has at most one degenerate eigenvalue, see \cite{SaSt}, Theorem 1). Therefore, Lemma 3 in \cite{SaSt} implies that there exist points $\zeta \in M_{T^{\theta}}$ arbitrarily close to $p$ with $H_{\rho}(X_{1}(\zeta),\overline{X_{1}(\zeta)}) > 0$. Because of the way $X_{1}$  was constructed, such a point $\zeta$ is a strictly pseudoconvex point. Consequently, for $\epsilon_{j} > 0$ and $\zeta$ near $p$, there exist real fields $Z_{\zeta,j} \in span_{\mathbb{R}}(ReX_{1},ImX_{1})$ of unit length so that $\mathcal{F}^{\epsilon_{j}}_{Z_{\zeta,j}}(z) \notin K$ for $z$ close enough to $\zeta$. This yields balls as in the theorem, but without control of the radii from below in terms of $\epsilon_{j}$. Because of the form of $Z_{\zeta,j}$, the uniform boundedness condition on the divergence of the fields $Z_{\zeta,j}$ is satisfied (near $p$). 

\medskip
We say that b$\Omega \setminus K$ satisfies a complex tangential cone condition if there exists a finite open real cone $\Gamma$ in $\mathbb{C}^n \approx \mathbb{R}^{2n}$ such that the following holds.  For each $p \in K$ there exists a complex tangential direction so that when $\Gamma$ is moved by a rigid motion to have vertex at $p$ and axis in that direction (that is, in the two dimensional real affine subspace determined by that direction), the (open) cone obtained intersects  b$\Omega$ in a set contained in b$\Omega \setminus K$. Theorem \ref{main} has the following corollaries.

\begin{corollary}\label{cone1}
Let $\Omega$ be a $C^\infty$-smooth bounded pseudoconvex domain in $\mathbb{C}^n$. Denote by $K$ the set of boundary points of infinite type.
For all points $\zeta$ in a neighborhood of $K$ in $b\Omega$, denote by $\lambda_{0}(\zeta)$ the smallest eigenvalue of the Levi form. 
Assume that b$\Omega$ satisfies the following conditions.  
There exist smooth complex tangential vector fields
$X_1, \dots, X_m$, defined on $b\Omega$ near K so that
$H_{\rho}(X_i(\zeta), \overline{X_i(\zeta)}) \leq C ~ {\lambda}_0(\zeta)$,
for some constant $C$ and all $\zeta$, such that $b\Omega \setminus K$ satisfies a complex tangential cone condition with 
the axis of the cone at $p \in K$ in $\text{span}_{\mathbb{R}}(\text{Re}X_1(p), 
\text{Im}X_1(p), \dots , \text{Re}X_m(p), \text{Im}X_m(p))$, for all $p \in K$.  
Then the $\bar{\partial}$-Neumann operator on $\Omega$ is compact.
\end{corollary}

In $\mathbb{C}^{2}$, the assumption in Corollary \ref{cone1} reduces to the simple requirement that $K$ satisfy a complex tangential cone condition; this is Corollary 2 in \cite{St2}. An example like the one described after the statement of Theorem \ref{main} shows that this is not sufficient in $\mathbb{C}^{n}$ when $n \geq 3$, not even when one assumes that the axis of the cone at $p \in K$ lies in the null space of the Levi form at $p$. So some complication in the statement of the corollary cannot be avoided. On the other hand, when the Levi form of b$\Omega$ has at most one degenerate eigenvalue, there is (locally) a complex tangential vector field $X_{1}$ and a constant such that $H_{\rho}(X_{1}(p), \overline{X_{1}}(p))=0$ and $H_{\rho}(X_{1}, \overline{X_{1}}) \leq C \lambda_0$ near $p \in K$ (see Example 3 above). With this additional information, it suffices to require that the axis of the cone lie in the null space of the Levi form.

\begin{corollary}\label{cone2}
Let $\Omega$ be a smooth bounded pseudoconvex domain in $\mathbb{C}^n$; assume that at each boundary point, the Levi form has at most one degenerate eigenvalue.If the set $b\Omega \setminus K$ satisfies a complex tangential cone condition with the axis of the cone at $p \in K$ in the null space of the Levi form at $p$, then the $\dbar$-Neumann operator on $\Omega$ is compact. 
\end{corollary}
 In $\mathbb{C}^{2}$, Corollary \ref{cone2} also reduces to Corollary 2 in \cite{St2}.

\section{Proof of Theorem \ref{main}}

We will establish a compactness estimate for forms in $C^{\infty}_{(0,1)}(\overline{\Omega}) \cap dom(\overline{\partial}^{*})$: for all $\epsilon > 0$, there is a constant $C_{\epsilon}$ such that $\|u\|^{2} \leq \epsilon(\|\overline{\partial}u\|^{2} + \|\overline{\partial}^{*}u\|^{2}) + C_{\epsilon}\|u\|_{-1}^{2}$. This is equivalent to compactness of $N_{1}$ (see for example \cite{FS1}, Lemma 1.1).

First note that we may assume that the vector fields $X_{k}, 1 \leq k \leq m$, are defined on all of $b\Omega$, by multiplying them with suitable cutoff functions that are identically equal to one near $K$. This preserves (\ref{new-cond}) (i.e. (\ref{new-cond}) now holds on all of $b\Omega$). Then we can extend the vector fields $Z_{p,j}$ and $X_{k}$ from b$\Omega$ to the inside of $\Omega$ by letting them be constant along the real normal, so that the extended fields, still denoted by $Z_{p,j}$ and $X_{k}$, respectively, are complex tangential to the level sets of the boundary distance. Finally, by multiplying by a suitable cutoff function that equals one near the boundary, we may assume that the fields $X_{k}$ are defined and smooth on $\overline{\Omega}$.

There are two ideas in the proof. The first comes from \cite{St2} and says that near a point $p \in K$, the values of a form $u$ can be expressed by the values near $\mathcal{F}_{Z_{p,j}}^{\epsilon_{j}}(p)$ plus the integrals of $Z_{p,j}u$ along the integral curves of $Z_{p,j}$. If points near $\mathcal{F}_{Z_{p,j}}^{\epsilon_{j}}(p)$ are of finite type, the contribution from there can easily be estimated in the required manner by subelliptic estimates. Because the integrals of $Z_{p,j}u$ in the second contribution are over curves of length $\epsilon_{j}$ ($|Z_{p,j}| = 1$), a (small) factor $\epsilon_{j}$ appears (via the Cauchy-Schwarz inequality) when computing the relevant $\mathcal{L}^{2}$-norms. There are overlap issues, but these are handled by the uniformity built into the assumptions (for example the uniform bound on $divZ_{p,j}$). The second idea concerns control of $\|Z_{p,j}u\|^{2}$, or ultimately $\|X_{k}u\|^{2}$, by $\|\overline{\partial}u\|^{2} + \|\overline{\partial}^{*}u\|^{2}$. In $\mathbb{C}^{2}$, this can be done for any complex tangential field, via maximal estimates (see \cite{D1,St2}). In higher dimensions, additional assumptions are needed; that \eqref{new-cond} suffices is implicit in  work of Derridj (\cite{D1}).

The first part of the proof follows \cite{St2} verbatim, with only one obvious modification. Fix $\epsilon >0$ and choose $j$ big enough so that $\epsilon_{j} < \epsilon$. The arguments in \cite{St2}, pp.~705--708 give estimate \eqref{11} below, the only modification being the exponent of $\epsilon_{j}$, which depends on the dimension. This dimension dependence arises from a comparison of volumes argument used to resolve certain overlap issues; see the paragraph in \cite{St2} that starts at the bottom of page 707. Note that when $n=2$, $2n+2-2nC_{3}=6-4C_{3}$, as in \cite{St2}, equation (14). Combining (the analogues of) equations (6) and (14) in \cite{St2} gives
\begin{equation}\label{11}
\begin{split}
\int\limits_{\Omega} \lvert u \rvert ^2  & \leq 2 \epsilon (\lVert \dbar u \rVert^2_0                                            + \lVert \dbars u \rVert^2_0 ) +     C_{\epsilon} \lVert u \rVert^2_{-1} 
 + 4 \epsilon_j \int_0^{\epsilon_j} \bigg( C(C_2)(\epsilon_j)^{2n-2nC_3}
2m \times \\  
                                         &\ \ \ \ \ \ \ \ \ \ \ \ \ \ \ \ \ \ \ \ \ \ \ \ \ \ \ \ \ \ \qquad \sum_{k=1}^{m} \underset{\Omega}{\int} \bigg[ 
\lvert \text{Re}X_k u(y) \rvert^2 + \lvert \text{Im}X_k u(y) \rvert^2\bigg] 
dV(y)\bigg) dt \\
                                         &\leq 2 \epsilon (\lVert \dbar u \rVert^2_0 + \lVert \dbars u \rVert^2_0 ) + 
C_{\epsilon} \lVert u \rVert^2_{-1} 
+ 16mC(C_{2})(\epsilon_{j})^{2n+2-2nC_{3}} \;\times \\
                                         &\ \ \ \ \ \ \ \ \ \ \ \ \ \ \ \ \ \ \ \ \ \ \ \ \ \ \ \ \ \ \ \ \ \ \qquad  \sum_{k=1}^{m} \underset{\Omega}{\int}  \bigg( \lvert X_k u(y)\rvert^2 +  \lvert \overline{X_k} u(y)\rvert^2\bigg) dV(y)\bigg) dt \;. \\
\end{split}
\end{equation}

In $\mathbb{C}^2$, the last term on the right hand side of \eqref{11} can be estimated using maximal estimates (\cite{D1}):

\begin{equation}\label{max}
\lVert \overline{X_k} u \rVert^2 + \lVert X_k u \rVert^2
\leq C_{k} \left( \lVert \dbar u \rVert_0^2 + \lVert \dbars u \rVert_0^2 \right).
\end{equation}

We are going to show that \eqref{max} also holds under hypothesis \eqref{new-cond} on the vector fields $X_{k}$. The estimates below owe much to \cite{D1}. For notational convenience, fix $k$ and denote $X_k$ by $X$. First note that the estimate on $\|\overline{X}u\|^{2}$ follows directly from the Kohn-Morrey formula (\cite{CS}, Proposition 4.3.1 with weight $e^{-\phi} \equiv 1$).
For $\lVert Xu \rVert^2$, we integrate by parts to obtain  
\begin{equation}\label{14}
\begin{split}
      \lVert Xu \rVert^2 &= - \int_{\Omega} u X \overline{Xu} 
                           + O(\|u\|\thinspace\|Xu\|)\\
                         &= -\int_{\Omega} u \left[X,\overline{X} \right]
                               \overline{u} - \int_{\Omega}u\overline{X}X\overline{u}  + O( \lVert u \rVert                            \thinspace \lVert Xu \rVert )\\
                         &= -\int_{\Omega} u\left[X,\overline{X} \right] 
                            \overline{u} + \lVert \bar{X}u \rVert ^2 +  O(\|u\|\thinspace\|Xu\| + \|u\|\thinspace\|\overline{X}u\|) \;.\\
\end{split}
\end{equation}
Using that $\|u\|^{2} \lesssim(\|\overline{\partial}u\|^{2} +  \|\overline{\partial}^{*}u\|^{2})$, $|\|u\|\thinspace\|Xu\| \leq l.c\|u\|^{2} + s.c\|Xu\|^{2}$ (\cite{CS}, (4.4.6) on page 79), and again the Kohn-Morrey formula for $\|\overline{X}u\|^{2}$, we see that it suffices to estimate the first term on the last line of \eqref{14}. Here, $s.c$ and $l.c.$ denote small and large constants, respectively; the term $\|Xu\|^{2}$ on the right hand side in \eqref{14} can then be absorbed into the left hand side.

Near $b\Omega$ (we may assume on the support of the cutoff functions used at the beginning of this section to extend $X$ to all of $\Omega$), we write $\left[X,\overline{X}\right ]\overline{u}$ as $H_{\rho}(X,\overline{X})(L_{n}-\overline{L_{n}})\overline{u} + A\overline{u} + \overline{B}\overline{u}$, where $A$ and $B$ are smooth complex tangential fields and $L_{n}$ is the complex normal $\sum_{j=1}^{n}(\partial\rho/\partial \overline{z_{j}})\partial / \partial z_{j}$ appropriately normalized. The contributions coming from $A\overline{u}$ and $L_{n}\overline{u}$, respectively, (or $\overline{A}u$ and $\overline{L_{n}}u$) are estimated as above. In the contribution from $\overline{B}\overline{u}$, we integrate $\overline{B}$ by parts, and proceed as before. We are left with estimating the term involving $H_{\rho}(X,\overline{X})\overline{L_{n}}\overline{u}$. We integrate $\overline{L_{n}}$ by parts (with a boundary term) to obtain
\begin{equation}\label{15}
\left |\int_{\Omega}uH_{\rho}(X,\overline{X})\overline{L_{n}}\overline{u} \right | \lesssim
\|\overline{L_{n}}u\|\thinspace\|u\| + \int_{b\Omega}H_{\rho}(X,\overline{X})|u|^{2}
+ O(\|u\|^{2}) \;.
\end{equation}
The first and the third term on the right hand side of \eqref{15} are estimated as above, while for the middle term we have
\begin{equation}\label{16}
\int_{b\Omega}H_{\rho}(X,\overline{X})|u|^{2} \lesssim \int_{b\Omega}\lambda_{0}|u|^{2} \lesssim \int_{b\Omega}H_{\rho}(u,\overline{u}) \lesssim \|\overline{\partial}u\|^{2} + \|\overline{\partial}^{*}u\|^{2} \; .
\end{equation}
The first inequality in \eqref{16} comes from \eqref{new-cond}, the last one from the Kohn-Morrey formula. We have slightly abused notation in the third term in \eqref{16}: $u$ is a $(0,1)$-form, not a vector field, but it is identified with a vector field in the usual way via its coefficients. Combining \eqref{14}, \eqref{15}, and \eqref{16} shows that \eqref{max} holds under the assumptions of Theorem \ref{main} (i.e. under \eqref{new-cond}). Inserting \eqref{max} into \eqref{11} gives
\begin{equation}\label{17}
\int_{\Omega}|u|^{2} \leq C(\epsilon + \epsilon^{2n+2-2nC_{3}})\left (\|\overline{\partial}u\|^{2} + \|\overline{\partial}^{*}u\|^{2}\right) + C_{\epsilon}\|u\|_{-1}^{2} \;,
\end{equation}
with $C$ independent of $\epsilon$. Since $2n+2-2nC_{3} > 0$ (because $C_{3} < (n+1)/n$ ), $\lim_{\epsilon \rightarrow 0^{+}}C(\epsilon + \epsilon^{2n+2-2nC_{3}}) = 0$, so \eqref{17} implies the compactness estimate we set out to prove. This completes the proof of Theorem \ref{main}.

\bigskip

\providecommand{\bysame}{\leavevmode\hbox to3em{\hrulefill}\thinspace}

\end{document}